\documentclass{article}
\usepackage[utf8]{inputenc}
\usepackage[left=3cm, right=3cm, top=3cm]{geometry}
\usepackage{setspace}
\usepackage{amsmath}
\usepackage{cool}
\usepackage{amssymb}
\usepackage{graphicx}
\usepackage{color}
\usepackage{mathrsfs}
\usepackage{enumitem}

\usepackage{algorithm}
\usepackage{mathtools}
\usepackage[noend]{algpseudocode}
\makeatletter
\def\BState{\State\hskip-\ALG@thistlm}
\makeatother
\usepackage{csvsimple}
\usepackage[short]{optidef}
\usepackage{pgfplotstable,filecontents}
\pgfplotsset{compat=1.9}
\usepackage[toc,page]{appendix}

\title{A Proof of Principle: Multi-Modality Radiotherapy Optimization}
\author{Roman Levin, Aleksandr Y. Aravkin, Minsun Kim}
\date{November 2019}

\begin{document}

\maketitle

\begin{abstract}
    Radiotherapy is used to treat cancer patients by damaging DNA of tumor cells using ionizing radiation. Photons are the most widely used radiation type for therapy, 
    having been put into use soon after the first discovery of X-rays in 1895.  
    However, there are emerging interests and developments of other radiation modalities such as protons and carbon ions, owing to their unique biological and physical characteristics that distinguish these modalities from photons. 
    Current attempts to determine an optimal radiation modality or an optimal combination of multiple modalities are empirical and in the early stage of development. 
    In this paper, we propose a mathematical framework to optimize full radiation dose distributions and fractionation schedules of multiple radiation modalities, 
    aiming to maximize the damage to the tumor while limiting the damage to the normal tissue to the corresponding tolerance level.  This formulation gives rise to a non-convex, mixed integer program and we propose a bilevel optimization algorithm, to efficiently solve it. The upper level problem is to optimize the fractionation schedule using the dose distribution optimized in the lower level. We demonstrate the feasibility of our novel framework and algorithms in a simple 2-dimensional phantom with two different radiation modalities, where clinical intuition can be easily drawn.  The results of our numerical simulations agree with the clinical intuition, validating our approach and showing the promise of the framework for further clinical investigation.
\end{abstract}

\section{Introduction}

The number of patients diagnosed with cancer is increasing every year, projecting over 1.9 million new cancer cases in 2020 according to the Centers for Disease and Control and Prevention (CDC) \cite{weir2015past}. More than half of all cancer patients go through radiotherapy in the course of their cancer treatment.  Radiotherapy utilizes 
\textcolor{black}{ionizing} radiation to kill cancer cells and is used as the primary treatment modality for certain cases or as an (neo)adjuvant modality before or after other treatment modalities such as surgery and chemotherapy. Although radiation kills cancer cells, it also damages normal tissue that is on its path.  Therefore, the goal of radiotherapy 
\textcolor{black}{is to maximize} the differential in the damages between the tumor and normal tissue.

\begin{figure}
\begin{centering}
\begin{tabular}{cc}
\includegraphics[width=3in]{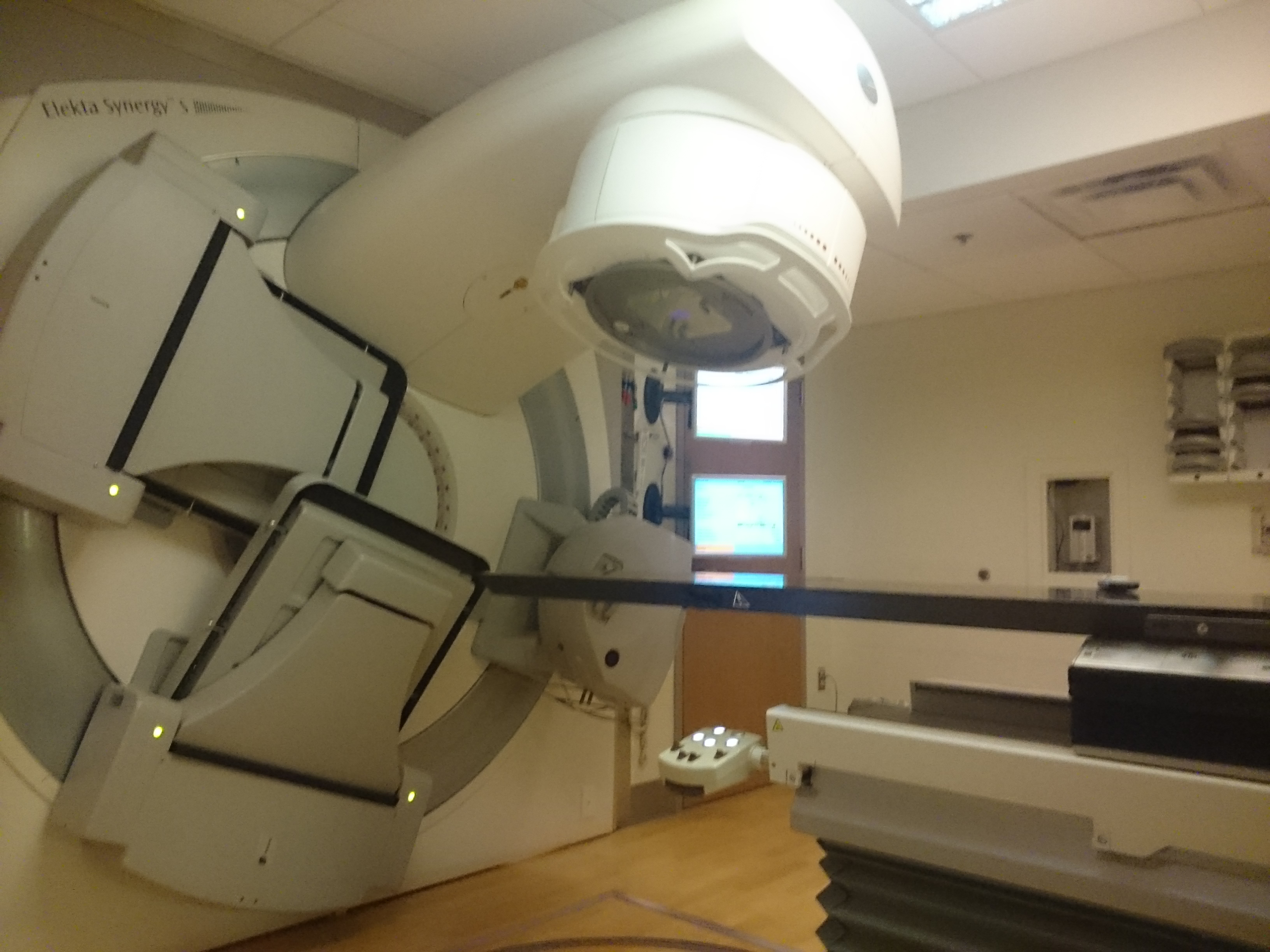}&\includegraphics[width=3.2in]{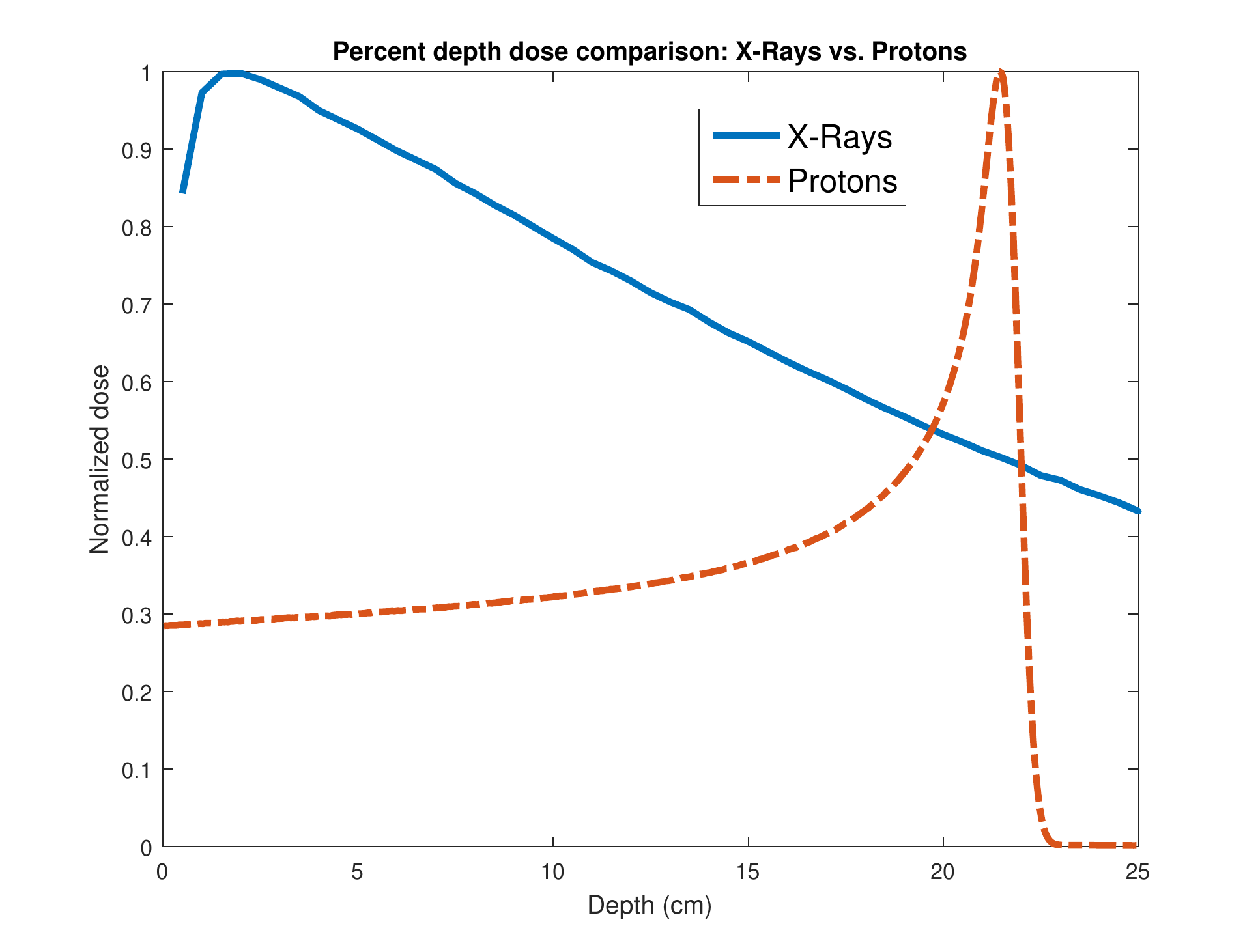}
\end{tabular}
\caption{(a) Elekta medical linear accelerator at the University of Washington Medical Center (b) Percent depth dose comparison between X-rays and protons..}
\label{fig:linac}
\end{centering}
\end{figure}

\begin{figure}
\begin{centering}
\begin{tabular}{c}
\includegraphics[width=5in]{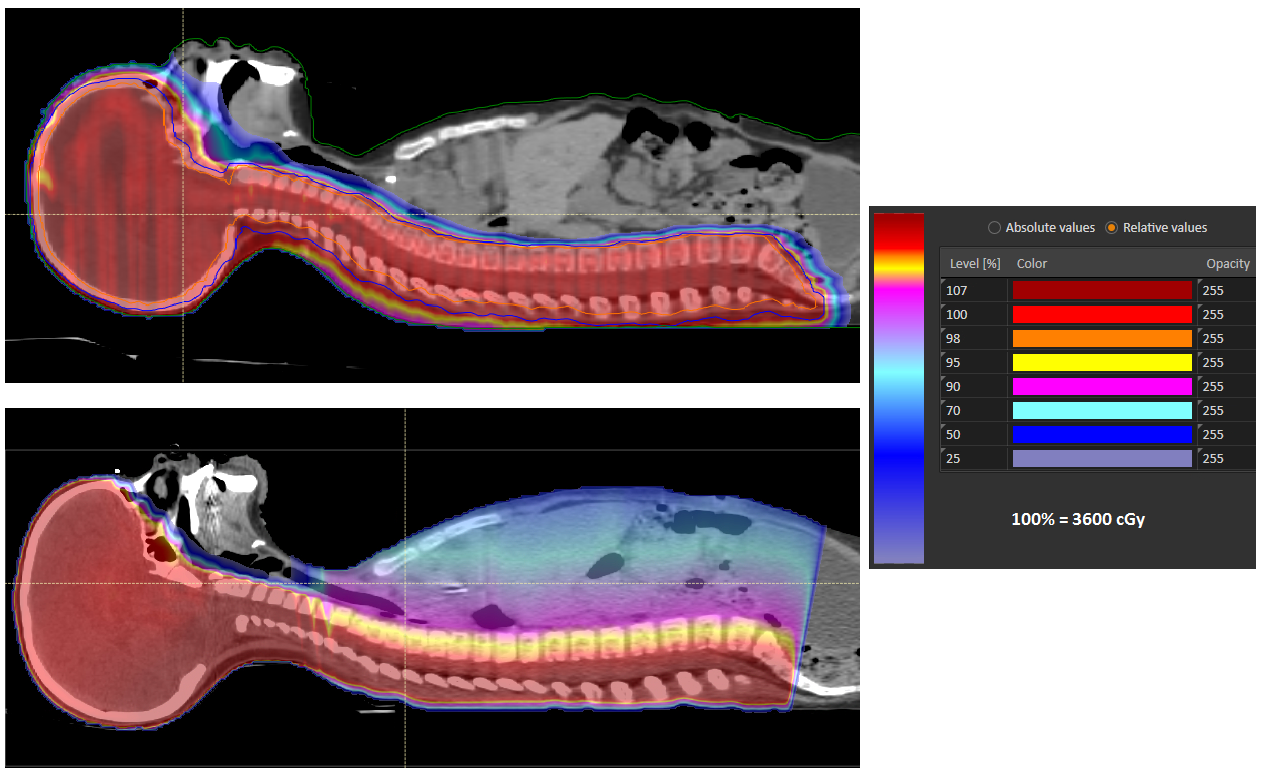}
\end{tabular}
\caption{Craniospinal irradiation patient plans (Top) Isodose distributions with a single posterior-anterior photon beam (Bottom) Isodose distributions with a single posterior-anterior proton beam.  Both plans are normalized to 3600 cGy at 100 \%.}
\label{fig:isodose}
\end{centering}
\end{figure}

External beam radiation therapy (EBRT) is a non-invasive type of radiotherapy, where the radiation generated by linear accelerators or cyclotrons is targeted at the tumor from outside the patient's body. Figure \ref{fig:linac} (a) shows a linear accelerator that produces X-rays and electrons to treat the patient lying on the table. There are currently multiple radiation types used in EBRT. The most widely used radiation type is photons (X-rays) and the therapeutic effect of photons on the tumor and normal tissue is well-established due to their long history of use.  However, there are emerging interests in the EBRT using heavy charged particles such as protons and carbon ions due to their unique dosimetric and biological characteristics. For example, Figure \ref{fig:linac} (b) shows the percent depth dose (PDD) of photons (X-rays) and protons.  Unlike photons, which deposit the maximum dose near the patient's surface, typically within 3 cm from the beam entrance for photons with the energy of less than 18 MV used in clinical practice, protons deposit the maximum dose at the end of
\textcolor{black}{their} range.  Since the range of protons is dependent on 
\textcolor{black}{their} energy, we can theoretically aim to deposit a large amount of dose in the tumor, leaving the normal tissue behind it almost no dose by adjusting the energy of the particles. This superior dosimetric effect of heavy charged particles compared to conventional photons can be easily seen in the isodose distributions in \textcolor{black}{Figure \ref{fig:isodose}} for the patient with craniospinal irradiation. The patient treated with protons on the top receives much less dose to the normal tissue outside the target area (whole brain and spinal cord) compared to the patient on the bottom, where the patient is treated with photons.  It is particularly beneficial for pediatric patients whose secondary malignancy is a great concern.

On the other hand, certain radiation types such as neutrons and carbon ions have superior biological effects compared to photons.  A relative biological effectiveness (RBE) is defined as the ratio of the absorbed dose of a reference radiation type (often photons) and another radiation type, which kills the same number of cells.  For example, the RBE of neutrons is defined as $D_{\text{X}}/D_{\text{n}}$, where $D_{\text{X}}$ and $D_{\text{n}}$ are the dose that kills the same number of cells using photons and neutrons respectively.  Therefore, a radiation type with a higher RBE kills more cells than photons given the same physical dose (absorbed energy per unit mass).  By definition, RBE is dependent on the radiation type, tissue type (what type of cells is considered in computing RBE), and the environment condition such as the level of oxygen \cite{Hall2006}. For example, fast neutrons used for radiotherapy have an RBE of 2-5 \cite{Britten2001,Wambersie1999} and it is particularly high in a hypoxic condition such as in the tumor. It is noteworthy that radiation with a higher RBE implies that it damages both tumor and normal tissue more than conventional photons but the differential is not uniformly scaled throughout various tissue types. {Therefore, there is an opportunity to exploit the differential to maximize the therapeutic effect.}  The summary of radiation types currently used in practice worldwide and their dosimetric and biological effect relative to conventional photons is presented in Table \ref{tab:radiationtypes}. 

\begin{table}[h!]
    \centering
    \begin{tabular}{c c c c}
    \hline
       Radiation Type  & Dosimetric effect & Biological damage & Cost \cite{Peeters2010} \\
       \hline
        Protons & Superior & Similar  & 3-4 times more \\
        Neutrons  & Similar & Superior & Similar \\
        Carbon ions & Superior & Superior &  5-6 times more\\
        \hline
    \end{tabular}
    \caption{Comparison of various radiation types. All are relative to photon (X-ray) external beam radiotherapy, which is currently the most widely used radiation type in practice.}
    \label{tab:radiationtypes}
\end{table}

Due to the unique dosimetric and biological effect of each radiation type, and practical considerations, there is no single modality that is superior to others in all aspects. Furthermore, it is not obvious which radiation type is optimal for a specific patient. Current efforts to determine an optimal radiation type or combination of radiation types are mostly empirical and largely depend on clinical intuitions \cite{Zhou2018,Haefner2018,Chi2017,Vogel2017}. More recently, there are some efforts to systematically optimize proton and photon treatments combined. We previously studied a simplified scenario, where the dose from two different radiation types is expressed in a scalar form, and proposed a mathematical framework to find an optimal fractionation for each radiation modality \cite{Nourollahi2018a} and its robust counter part in \cite{Nourollahi2018b}. Unkelbach et al. minimized the mean biologically effective dose (BED) of organs-at-risk while prescribing a fixed BED to the tumor with a fixed fractionation for both protons and photons \cite{Unkelbach2018}. Gao et al. studied the hybrid proton-photon inverse planning optimization in \cite{Gao2019}, where they also optimized the dose distribution from both modalities using the fixed prescription dose to the tumor and fixed fractionation schedules.  Eikelder et al. studied the fluence map-fractionation schedule optimization of proton and photon combined using a sparing factor \cite{TenEikelder2019}.  The idea of a sparing factor in the fractionation schedule optimization is explained in detail in \cite{Saberian2015} with a single modality case. In summary, the relative dose of OAR in each voxel is fixed as a fraction of the tumor dose through a sparing factor. Eikelder et al. took a heuristic approach to proton-photon modality fractionation optimization problem, where the maximum feasible BED to the tumor was computed for all possible combinations of fractions for each modality to find an optimal fractionation schedule. In particular, their approach separates the fractionation schedule optimization from the fluence map optimization because the sparing factors are not optimization variables.

\textcolor{black}{The purpose of this study is to:
\begin{enumerate}
    \item Set up a full-scale, rigorous mathematical framework to simultaneously optimize dose distributions and fractionations of two or more radiation modalities.  In other words, the dose distribution optimization and fractionation optimization of multiple modalities are integrated.
    \item Develop an efficient optimization algorithm to find optimal fractionation and corresponding optimal fluence maps (and resulting dose distributions) for each modality.
    \item Test the feasibility and clinical relevance of the proposed  framework on a small-scale phantom, where clinical intuition can be used to validate the results of the numerical simulation.
\end{enumerate}}

\subsection{Contribution}

We develop a rigorous optimization approach for the radiation treatment planning problem using $M$ radiation modalities.  We propose a novel, integrated framework that can simultaneously optimize dose distributions and fractionations of multiple modalities in full, voxel-by-voxel scale. In other words, we seek the solutions to the problem, which exploit the full flexibility of treatment parameters (i.e. number of modalities, fluence maps, number of fractions).

Our work focuses on developing a non-convex optimization framework, where the biological effect (BE) to the tumor is maximized while the normal tissue BE is kept under the tolerance level. The optimization variables are the fluence maps of each modality ($u_m, m=1, \cdots, M$) and its corresponding number of fractions ($N_m$). We also develop a bilevel optimization algorithm to solve for optimal $\{u_m\}_{m=1}^M$ and $\{N_m\}_{m=1}^M$.  The lower level optimizes $u_m$ using non-convex relaxation for a fixed fractionation schedule. The upper level relaxes the integrality of $N_m$ and solves the resulting continuous optimization problem using a trust region algorithm, calling the lower level as a subroutine. We demonstrate the feasibility of the proposed framework by applying it to a small-scale phantom case. 

This paper is organized as follows: the problem formulation and optimization framework are presented in Section \ref{sec:prob} followed by the optimization algorithms to solve the problem in Section \ref{sec:solution}. The numerical simulations and their results are presented in Section \ref{sec:numerical}.  We conclude in Section \ref{sec:conclusion} with some remarks for the future work.

\section{Problem Formulation}
\label{sec:prob}

Let us consider optimizing the dose distributions using $M$ modalities, where each modality $m$ delivers $N_m$ fractions ($m=1,2, \cdots, M$).  Biological effect (BE) based on the linear-quadratic (LQ) cell-survival model is  widely used in radiation oncology to characterize the effect of the physical dose (the energy absorbed per unit mass) combined with the fractionation effect \cite{Hall2006}.  We 
use BE to compare the effect of two different fractionation schedules on cell-killing since the same total (physical) dose could lead to a different biological outcome depending on the fractionation schedule.  For example, 10 Gy delivered to the tumor in 1 fraction kills more tumor cells than 10 Gy delivered in 5 fractions. 
\textcolor{black}{BE of delivering one-dimensional dose $d_m$ to the tumor for $N_m$ fractions for a single radiation modality $m$ is given by}
\begin{eqnarray}
\label{eqn:BE}
\mbox{BE} = \alpha_m N_m d + \beta_m N_m d^2 - \tau(N_m),
\end{eqnarray}
where $\alpha_m$ and $\beta_m$ are radiation modality and tissue-specific radiobiological parameters in the LQ model and $\tau(N_m)$ is a tumor proliferation term, only depends on the length of treatment.

The goal is to maximize the total tumor BE
(sum of BEs across modalities)
while keeping each of the organ-at-risk (OAR) total BEs under tolerance, using $M$ radiation modalities with $N_m$ fractions for each modality. In the multidimensional case, the radiobiological parameters of each radiation modality $\alpha_m$ and $\beta_m$ are  vectors, and each dose distribution $d_m$ is the image of fluence map $u_m$ under the dose mapping. We use a linear dose mapping for the dose calculation, so that $d_m = A_m u_m$.  OAR constraint types included in this paper are mean dose constraints and maximum dose constraints, which are common in practice for parallel\footnote{The organ remains functional when part of it is damaged by radiation.} and serial\footnote{If any part of an organ is damaged, the organ becomes no longer functional.} normal tissue types respectively. Another common constraint type for OAR in practice is dose-volume (DV) constraints, which specify a critical volume of OAR that must receive less than a certain critical dose value. DV constraints can be easily handled using a constraint generation method, where the initial optimization is done without DV constraints and then the maximum constraints are applied to specific voxels in the subsequent optimization if DV constraints are violated in the first optimization \cite{Kim2015, Saberian2015, Saberian2017}. Since DV constraints can be handled by maximum constraints with some modification in the algorithm, we do not include DV constraints for the brevity of notations.

\subsection{Notation}
We use the following notation to describe the proposed model:
\begin{itemize}
    \item $M$: total number of radiation modalities 
    \item \textcolor{black}{$u_m$: fluence map (beamlet intensities) for modality $m$}
    \item $N_m$: number of fractions for modality $m$
    \item $T_m$: tumor dose coefficient matrix for modality $m$ such that $T_m u_m$ 
    \textcolor{black}{gives the dose distribution} delivered to the tumor using a fluence map $u_m$
    \item $T_m(j)$: $j$-th row of matrix $T_m$
    \item $H_m^i$: $i$-th OAR dose coefficient matrix for modality $m$
    \item $H_m^i(j)$: $j$-th row of matrix $H_m^i$
    \item $n$: number of OARs
    \item $l$: number of voxels in the tumor 
    \item $J_0$: \textcolor{black}{voxel index set for the tumor (so $|J_0| = l$)}
    \item $J_i$: \textcolor{black}{voxel index set for the $i$-th OAR, $i=1,2,\cdots, n$}
    \item $I_{\text{mean}}$: index set \textcolor{black}{of} OARs with mean-dose constraint
    \item $I_{\text{max}}$: index set \textcolor{black}{of} OARs with max-dose constraint
    \item $\alpha_m$, $\beta_m$: vectors of the linear and quadratic \textcolor{black}{radiobiological} coefficients in the tumor BE for modality $m; \alpha_m, \beta_m\in \mathbb{R}_{++}^{|J_0|}$
    \item $\gamma_m^i$, $\delta_m^i$: vectors of the linear and quadratic \textcolor{black}{radiobiological} coefficients in the $i$-th OAR BE for modality $m; \gamma_m, \delta_m\in \mathbb{R}_{++}^{|J_i|}$
    \item $\alpha_m(j)$, $\beta_m(j)$, $\gamma_m^i(j)$, $\delta_m^i(j)$: $j$-th element of the corresponding vector
    \item $C^i_{\text{mean}}$: tolerance BE for the $i$-th OAR with a mean dose constraint
    \item $C^i_{\text{max}}$: tolerance BE for the $i$-th OAR with a maximum dose constraint
\end{itemize}

\subsection{Optimization Framework}
\textcolor{black}{We now formulate the fully general multi-modality radiotherapy framework:
$M$ radiation modalities with $N_m$ fractions for the $m$-th modality. We denote tumor dose coefficient matrices  which map beamlet intensities $u_m$ to the dose distribution for each modality by $T_m$, and OAR dose coefficient matrices that map $u_m$ to the dose distribution delivered to the $i$-th OAR by $H^i_m$. Our objective is to maximize the total tumor BE, subject to the mean and max constraints on the OAR BEs:}

\textcolor{black}{
\begin{maxi!}[1]
  {\{u_m, N_m\}_{m=1}^M}{\sum_{m=1}^M \sum_{j \in J_0} N_m \alpha_{m}(j)(T_m(j) u_m) + N_m \beta_m(j)(T_m(j)u_m)^2- \tau(N_m)}{}{(\text{P0})\quad}\tag{$\text{P0}.1$} \label{eqn:gen-framework} 
  \addConstraint{\sum_{m=1}^M \sum_{j \in J_i} N_m \gamma^i_{m}(j)(H^i_m(j) u_m) + N_m \delta^i_m(j)(H^i_m(j)u_m)^2}{\le C^i_{\text{mean}}, \quad}{\forall i \in I_{\text{mean}}}\tag{$\text{P0}.2$}\label{eqn:meandoseconstraint}
  \addConstraint{\max_{j \in J_i} \left\{\sum_{m=1}^M N_m \gamma^i_{m}(j)(H^i_m(j) u_m) + N_m \delta^i_m(j)(H^i_m(j)u_m)^2\right\}}{\le C^i_{\text{max}}, \quad}{\forall i \in I_{\text{max}}}\tag{$\text{P0}.3$}\label{eqn:maxdoseconstraint}
  \addConstraint{1 \le \sum_{m=1}^M N_m}{\le N_{\text{max}}}{}\tag{$\text{P0}.4$}
  \addConstraint{N_m}{\in \mathbb{Z}_{\ge 0}, \quad}{m = 1, \ldots, M} \tag{$\text{P0}.5$}
  \addConstraint{u_m}{\succeq 0, \quad}{m = 1, \ldots, M}\tag{$\text{P0}.6$}
\end{maxi!}
}

This problem is a non-convex mixed integer program, which is NP-hard in general. We tackle this problem by devising a bilevel optimization algorithm, where in the upper level, we relax integrality constraints and optimize the fractionation schedule over $\{N_m\}_{m=1}^M$  using the optimal fluence map $\{u_m\}_{m=1}^M$ obtained from the lower level for a given fractionation schedule. The details of the  algorithms are presented in the following section.

\section{Optimization Algorithms}\label{sec:solution}
In this section we develop a bilevel optimization algorithm to solve Problem (P0). The upper level optimizing over fractionation schedules ($\{N_m\}_{m=1}^M$) is detailed in Section \ref{subsec:upper} and the lower level approach (used as a subroutine in the upper level) to compute the optimal fluence map ($\{u^*_m\}_{m=1}^M$) for fixed \textcolor{black}{$\{N_m\}_{m=1}^M$} is presented in Section \ref{subsec:lower}.

\subsection{Upper Level: Fractionation Schedule Optimization}
\label{subsec:upper}

In the upper level, we optimize the number of fractions of each modality. We first convert the maximization problem in Problem (P0) to a minimization problem.  The tumor proliferation term in Equations (\ref{eqn:BE}) and (P0.1) depends on the total treatment length, $N$. Assuming that there is no tumor lagging time\footnote{time it takes for the tumor to start proliferation after the treatment starts}, the tumor proliferation term with the total $N$ fractions can be defined as
\begin{eqnarray}
\label{eqn:prolif}
\tau(N) = \frac{\ln 2 (N-1)}{T_d},
\end{eqnarray}
where $T_d$ is the tumor doubling time\footnote{time it takes for the tumor to double in the number of cells} \cite{Hall2006}. Using~\eqref{eqn:prolif}, we can rewrite ($\text{P0}$) as follows:
\\
\textcolor{black}{
\begin{mini!}[1]
  {\{u_m, N_m\}_{m=1}^M}{\sum_{m=1}^M N_m\left\{{-\alpha_m}^TT_iu_m - (T_mu_m)^T\mbox{diag}(\beta_m)(T_mu_m)\right\}
  +\frac{l\left(\sum_{m=1}^M N_m - 1\right)\ln{2}}{T_d}}{}{(\text{P1})\quad}\tag{$\text{P1}.1$} \label{eqn:upperlevel}
  \addConstraint{\sum_{m=1}^{M}N_m\left\{{\gamma_m^i}^T H_m^i u_m+(H_m^i u_m)^T \mbox{diag}(\delta^i_m) (H_m^i u_m)\right\}}{\le C^i_{\text{mean}}, \quad}{\forall i \in I_{\text{mean}}}\tag{$\text{P1}.2$}\label{eqn:meandoseconstraint2}
  \addConstraint{\max_{j \in J_i} \left\{\sum_{m=1}^M N_m \gamma^i_{m}(j)(H^i_m(j) u_m) + N_m \delta^i_m(j)(H^i_m(j)u_m)^2\right\}}{\le C^i_{\text{max}}, \quad}{\forall i \in I_{\text{max}}}\tag{$\text{P1}.3$}\label{eqn:maxdoseconstraint2}
  \addConstraint{1 \le \sum_{m=1}^M N_m}{\le N_{\text{max}}}{}\tag{$\text{P1}.4$}
  \addConstraint{N_m}{\in \mathbb{Z}_{\ge 0}, \quad}{m = 1, \ldots, M} \tag{$\text{P1}.5$}
  \addConstraint{u_m}{\succeq 0, \quad}{m = 1, \ldots, M}\tag{$\text{P1}.6$}
\end{mini!}
}
\noindent
where $l$ is the total number of voxels in the tumor.
Let $F(\{u_m\}_{m=1}^M, \{N_m\}_{m=1}^M) $ denote the objective function in (\ref{eqn:upperlevel}). 
To solve the problem, we relax integrality constraints and solve 
the continuous optimization problem. 

\begin{mini!}[1]
  {N_1,\cdots, N_M}{V(N_1,\cdots,N_M)}{}{(\text{P2})}\tag{P2.1}
  \addConstraint{1 \le \sum_{m=1}^M N_m}{\le N_{\text{max}}}{}\tag{P2.2}
   \addConstraint{N_m}{\ge 0, \quad}{m = 1, \ldots, M} \tag{P2.3},
\end{mini!}
    where $V(N_1, \dots, N_M)$ is the {\it value function} of $\{N_m\}_{m=1}^M$ defined by 
    \begin{equation}
        \label{eqn:valuefunc}
        V(N_1, \dots, N_M) := F(\{u^*_m(N_1, \dots, N_M)\}_{m=1}^M,N_1, \dots, N_M)
    \end{equation}
    where each $u^*_m(N_1, \dots, N_m)$ is the optimal fluence map solution 
    for fixed fractionation $(N_1, \dots, N_m)$:
    \begin{eqnarray}
    \label{eqn:lowerustar}
    \{u^*_m(N_1, \dots, N_M)\}_{m=1}^M &=& \arg\min_u (\text{P2} (N_1, \dots, N_M)).
    \end{eqnarray}
Every evaluation of the value function $V$ requires solving an optimization problem 
in the fluence maps (see Equation (\ref{eqn:lowerustar})), and this is done 
using the lower level soluton approach discussed in Section \ref{subsec:lower}. Problem (P2) has a 
nonlinear nonconvex objective with simple linear inequality constraints, and we solve it using a trust region method \cite{conn2000trust} for constrained optimization as implemented in Python package, SciPy \cite{scipy2019arXiv190710121V}. 
    The solution of ($\text{P2}$) is rounded to the nearest integers, $N^*_f$ in a post-processing step.  Once we have an integral $N^*_f$, we also 
    update fluence maps $u^*_f = \arg\min_u (\text{P1} (N^*_f))$ to ensure that the optimality and feasibility are enforced for the integer solutions. The upper level solution approach is described in Algorithm \ref{algo:optimize_N}. Since the problem is nonconvex, we repeat Algorithm~\ref{algo:optimize_N} for multiple initial guesses of $\{N^{(0)}_m\}$ and the best solution is chosen as the final optimal solution.
  
\begin{algorithm} 
\caption{$\{N_m\}$ Fractionation Schedule Optimization}\label{algo:optimize_N} 
\begin{algorithmic}[1]
\State \textbf{Input:} $u^{(0)}, N_1^{(0)}, \dots, N_M^{(0)}$
\Function{ObjectiveFun}{$u, N1, \dots, N_M$}
\State \Return $
 \sum_{i=1}^M N_i({\alpha_i}^TT_iu_i - (T_iu_i)^T\mbox{diag}(\beta_i)(T_iu_i)) + n_{tx}\left(\sum_{j=1}^M N_j - 1\right)\ln{2}/T_d.
$
\EndFunction
\Function{ValueFun}{$N_1, \dots, N_M$} \Comment{Define the value function to optimize}
\State $u^*_N \gets $\textsc{LowerLevelSolver}($u^{(0)}, N_1, \dots, N_M$)
\State \Return \textsc{ObjectiveFun}($u^*_N, N1, \dots, N_M$)
\EndFunction
\State  $N_1^{*}, \dots, N_M^{*} \gets $ \textsc{TrustRegionConstr}(\textsc{ValueFun}, $N_1^{(0)}, \dots, N_M^{(0)}$, $\sum_{j=1}^M N_j \le 25$, $\{N_j \ge 0\}_1^M$)
\State $u^* \gets$ \textsc{LowerLevelSolver}($u^{(0)}, [N_1^{*}], \dots, [N_M^{*}]$)
\State \textbf{Output:} $u^*, N_1^{*}, \dots, N_M^{*}$
\end{algorithmic}
\end{algorithm}


\subsection{Lower Level: Fluence Map Optimization for Fixed Fractionation}
\label{subsec:lower}
In this section, we describe the lower level solution required to compute~\eqref{eqn:lowerustar} 
for a given fractionation schedule. The tumor proliferation term is independent of $u_m$ 
and does not affect the problem for fixed $N_m$. 
Next, every maximum dose constraint in \eqref{eqn:maxdoseconstraint2} can be viewed as a mean dose constraint applied to every single voxel of a given OAR, essentially treating each 
of those voxels as a mean-dose OAR in its own right. The dose mapping matrices for those "new OARs" are comprised of the corresponding rows of the original OAR dose mapping matrices. Let us specify the dose mapping and BE coefficient matrices in terms of these new OARs. First, noting that all the voxel index sets 
$J_1, \dots, \textcolor{black}{J_{n}}$ are disjoint, define the index set of all OAR voxels with maximum dose constraints as follows: 
\begin{equation}
    \tilde I_{\text{max}} = \{i=(j,k)| k \in I_{\text{max}}, j \in J_k\}.
\end{equation}
$\tilde I_{\text{max}}$ is the index set of our new "mean-dose OARs". Now, for any OAR $i = (j,k) \in \tilde I_{\text{max}}$, define the corresponding dose coefficient matrix for modality $m$ as $H^k_m(j)$ (in fact, this would be a row-vector) and arrange $M$ dose coefficient matrices into a single block-diagonal generalized dose matrix $H^i$ for the OAR $i$. Stack the corresponding linear BE coefficients into vectors and include the fixed $N_m$ to obtain the following generalized linear BE coefficients:
\begin{align}
H^i = \begin{bmatrix}
H^k_1(j) & & & \\
& H^k_2(j) & & \\
&& \ddots &\\
&&&H^k_M(j)
\end{bmatrix}, 
\tilde{\gamma}^i = 
\begin{bmatrix}
N_1{\gamma}^k_1(j) \\
N_2{\gamma}^k_2(j) \\
\vdots \\
N_M{\gamma}^k_M(j)
\end{bmatrix}.
\end{align}
Similarly, for every OAR $i \in I_{\text{mean}}$ arrange the dose matrices of each modality into block-diagonal generalized dose matrices and stack the linear BE coefficients into large generalized linear BE coefficient vectors, do the same for the tumor dose matrices and BE coefficients getting: 
\begin{align}
H^i = \begin{bmatrix}
H^k_1 & & & \\
& H^k_2 & & \\
&& \ddots &\\
&&&H^k_M
\end{bmatrix}, 
\tilde{\gamma}^i = 
\begin{bmatrix}
N_1{\gamma}^k_1 \\
N_2{\gamma}^k_2 \\
\vdots \\
N_M{\gamma}^k_M
\end{bmatrix},
T = \begin{bmatrix}
T_1 & & & \\
& T_2 & & \\
&& \ddots &\\
&&&T_M
\end{bmatrix}, 
\tilde{\alpha} = \begin{bmatrix}
-N_1\alpha_1\\
-N_2\alpha_2\\
\vdots\\
-N_M\alpha_M
\end{bmatrix}.
\end{align}
{Rearranging the sums in equations \eqref{eqn:gen-framework} and  \eqref{eqn:meandoseconstraint} of the Problem (P0) and rewriting the max-dose constraints \eqref{eqn:maxdoseconstraint} in terms of the new one-voxel mean-dose-constrained OARs}, we reformulate the Problem (P0) for fixed  \textcolor{black}{$\{N_m\}_{m=1}^M$} as follows: 
\textcolor{black}{
\begin{mini!}[1]
  {u}{{\tilde{\alpha}}^T(Tu) - (Tu)^T B(Tu)}{}{}\nonumber
  \addConstraint{\tilde{\gamma^i}^T H^i u+(H^iu)^TD^i(H^iu)}{\le C^i_{\text{mean}}, \quad}{\forall i\in I_{\text{mean}}}\nonumber
  \addConstraint{\tilde{\gamma^i}^TH^iu+(H^iu)^TD^i(H^iu)}{\le C^i_{\text{max}},\quad}{\forall i\in \tilde{I}_{\text{max}}}\nonumber
  \addConstraint{u}{\succeq 0 . \quad}{}\nonumber
\end{mini!}
}
where 
\textcolor{black}{
\[
 u = \begin{bmatrix}
u_1 \\
u_2 \\
\vdots \\
u_M
\end{bmatrix},
B = \begin{bmatrix}
N_1 \mbox{diag}(\beta_1) & & & \\
 & N_2 \mbox{diag}(\beta_2) & & \\
& & \ddots & \\
& & & N_M \mbox{diag}(\beta_M)\\
\end{bmatrix},
\]
and 
\[
D^i = \begin{bmatrix}
N_1 \mbox{diag}(\delta^i_1) & & & \\
 & N_2 \mbox{diag}(\delta^i_2) & & \\
& & \ddots & \\
& & & N_M \mbox{diag}(\delta^i_M)\\
\end{bmatrix}
\qquad \forall i\in I_{\text{mean}}, \]
\[
D^i = \begin{bmatrix}
N_1 {\delta}^k_{1}(j) & & & \\
 & N_2 {\delta}^k_{2}(j) & & \\
& & \ddots & \\
& & & N_M {\delta}^k_{M}(j)\\
\end{bmatrix}
\qquad \forall i = (j,k) \in \tilde I_{\text{max}}.\]}
Here, $\mbox{diag}(v)$ denotes the diagonal matrix formed from a vector $v$, and 
the matrix $B$ is block-diagonal with blocks $N_m\mbox{diag}(\beta_m)$. \textcolor{black}{Likewise, the matrices $D^i$ for $i\in I_{\text{mean}}$ are block-diagonal with blocks $N_m \mbox{diag}(\delta^i_m)$. Finally, $D^i$ for $i \in \tilde I_{\text{max}}$ corresponding to the one-voxel mean-dose-constrained OARs are diagonal $M \times M$ matrices.}
For brevity, we denote the quadratic form with the matrix $B$ as
\textcolor{black}{\[
f(x) := x^T B x.
\]}
We now restate the optimization formulation for the fixed $\{N_m\}_{m=1}^M$:

\textcolor{black}{
\begin{mini!}[1]
  {u}{{\tilde{\alpha}}^T(Tu) - f(Tu)}{}{(\text{P3})}\tag{P3.1}
  \addConstraint{\tilde{\gamma^i}^T H^i u+(H^iu)^TD^i(H^iu)}{\le C^i_{\text{mean}}, \quad}{\forall i\in I_{\text{mean}}}\tag{P3.2}
  \addConstraint{\tilde{\gamma^i}^TH^iu+(H^iu)^TD^i(H^iu)}{\le C^i_{\text{max}},\quad}{\forall i\in \tilde{I}_{\text{max}}}\tag{P3.3}
  \addConstraint{u}{\succeq 0.}{}\tag{P3.4}
\end{mini!} \label{eqn:general_fixed_N}}

\subsubsection{Non-Convex Relaxation}
\label{subsec:nonconvexrelax}
The optimization problem ($\text{P3}$) is non-convex since we are minimizing a concave objective function and it is difficult to solve directly. To attack this problem, 
we follow the ideas of~\cite{zheng2018relax} and introduce auxiliary variables, $w_0$ and $w_i$ with $i=1,2, \cdots, \tilde{n}$, which gives us a more tractable 
relaxed problem:
\textcolor{black}{
\begin{mini!}[1]
  {u, w_0, \{w_i\}_{i=1}^{\tilde{n}}}{\tilde{\alpha}^Tw_0 - f(w_0) + \frac{1}{2\eta_0} \|w_0 - Tu\|^2 + \sum_{i=1}^{\tilde{n}}\frac{1}{2\eta_i}\|w_i - H^iu\|^2}{}{(\text{P4})}\tag{P4.1}
  \addConstraint{\tilde{\gamma}^Tw_i+w_i^TD^iw_i}{\le C^i_{\text{mean}}, \quad}{\forall i\in I_{\text{mean}}}\tag{P4.2}
  \addConstraint{\tilde{\gamma}^Tw_i+w_i^TD^iw_i}{\le C^i_{\text{max}},\quad}{\forall i\in \tilde{I}_{\text{max}}}\tag{P4.3}
  \addConstraint{u}{\succeq 0.}{}\tag{P4.4}
\end{mini!}}
where $\tilde{n} = |I_\text{mean}|+|\tilde{I}_\text{max}|$.
The norm penalties force $w_0$ and $w_i$ to be close to $Tu$ and $H^iu$. The parameters $\eta_0$ and $\eta_i$ control the degree of closeness, and, as $\eta_0, \eta_1, \dots, \eta_{\tilde{n}}$ go to zero, we recover (P3) from (P4). 
We develop an automatic approach to select these parameters in Section \ref{sec:parameters}. By design, the auxiliary variables $w_i$ always meet the original mean or maximum BE constraints for every OAR. 

To solve ($\text{P4}$), we \textcolor{black}{use block-coordinate descent and} iteratively update $u, w_0$, and $w_i$s. We now describe each update in detail.

\subsubsection*{Update $w_0, w_i$}
For fixed $u$, the problem we solve is given by 

\textcolor{black}{
\begin{mini!}[1]
  {w_0, \{w_i\}_{i=1}^{\tilde{n}}}{\tilde{\alpha}^Tw_0 - f(w_0) + \frac{1}{2\eta_0} \|w_0 - Tu\|^2 + \sum_{i=1}^{\tilde{n}}\frac{1}{2\eta_i}\|w_i - H^iu\|^2}{}{(\text{P5})}\tag{P5.1}
  \addConstraint{\tilde{\gamma}^Tw_i+w_i^TD^iw_i}{\le C^i_{\text{mean}}, \quad}{\forall i\in I_{\text{mean}}}\tag{P5.2}
  \addConstraint{\tilde{\gamma}^Tw_i+w_i^TD^iw_i}{\le C^i_{\text{max}},\quad}{\forall i\in \tilde{I}_{\text{max}}}\tag{P5.3}
\end{mini!}}
This problem is decoupled in $w_0$ and $w_i$s and is therefore equivalent to solving for $w_0$ and $w_i$  separately. For $w_0$, dropping the constant terms, we have 
\begin{equation}
\min_{w_0} {\alpha}^Tw_0 - f(w_0) + \frac{1}{2\eta_0} \|w_0 - Tu\|^2,
\end{equation}
which is equivalent to
\begin{equation}
\label{prox}
\min_{w_0} -f(w_0) + \frac{1}{2\eta_0} \|w_0 - (Tu -\eta_0\tilde{\alpha})\|^2.
\end{equation}
The solution of this minimization step for $w_0$ can be written compactly as
\begin{eqnarray}
\label{eqn:w0update}
w_0^+ = \mbox{prox}_{-\eta_0f} (Tu - \eta_0\tilde{\alpha}),
\end{eqnarray}
where the {\it proximal operator}, or prox, is defined by
\begin{eqnarray}
\label{eqn:prox}
\mbox{prox}_{-\eta_0f}(y) = \arg\min_x \left\{-f(x) + \frac{1}{2\eta_0}\|x-y\|^2 \right\}.
\end{eqnarray}

Equation (\ref{eqn:w0update}) gives us the update $w_0^+$ for $w_0$. The proximal operator always exists for closed convex functions, but care must be taken in the nonconvex case, and in particular for the concave $-f(x)$. By analyzing the problem, we find the range of values $\eta_0$ for which 
the proximal operator is well-defined (see Appendix \ref{appendix:proximal}): 
\begin{eqnarray}
\label{eqn: prox_solution}
(\mbox{prox}_{-\eta_0f}(y))_j = \begin{cases} 
\infty, & \max_i(B_{ii}-\frac{1}{2\eta_0}) > 0 \;\&\; j = \arg\max_i(B_{ii}-\frac{1}{2\eta_0})\\
0, & \max_i(B_{ii}-\frac{1}{2\eta_0}) > 0 \;\&\; j \not= \arg\max_i(B_{ii}-\frac{1}{2\eta_0})\\
\frac{1}{\eta_0}y_j/(-2B_{jj}+\frac{1}{\eta_0}), & \max_i(B_{ii}-\frac{1}{2\eta_0}) \le 0.
\end{cases}
\end{eqnarray}
In order for $\mbox{prox}_{-\eta_0 f}$ to be well-defined, we must have $\max_i(B_{ii}-\frac{1}{2\eta_0}) \le 0$. This assumption forces a lower limit on 
the penalty parameter, making sure problems (P3) and  (P4) are close, and gives us a starting point for the parameter selection process discussed in Section \ref{sec:parameters}.

Next, we consider the optimization problem with respect to each $w_i$ keeping only those terms that depend on $w_i$. 
We have $\tilde{n}$ number of problems with the same structure:
\begin{mini!}[1]
  {w_i}{\|w_i - H^iu\|^2}{}{(\text{P6})}\tag{P6.1}
  \addConstraint{\tilde{\gamma}^Tw_i+w_i^TD^iw_i}{\le C^i_{\text{mean/max}}, \quad}{\forall i\in I_{\text{mean}}/\tilde I_{\text{max}}}\tag{P6.2}
\end{mini!}


The solution of this problem is the projection of $H^iu$ onto the convex set $\Omega_i = \{w_i: \tilde \gamma^Tw_i+w_i^TDw_i \le C^i_{\text{mean/max}}\}$. The solution method to find the projection,
\begin{eqnarray*}
&\mbox{proj}_{\Omega_i}(v) = \arg\min_{w \in \Omega_i} \|w - v\|^2 
\end{eqnarray*}
is presented in Appendix \ref{appendix:projection}. Projection onto a closed convex set is a special 
case of the prox operator, and is always well-defined and single-valued. 

\subsubsection*{Update $u$}
We now consider the subproblem for $u$ for fixed $w_0$ and $w_i$s. Dropping the constant terms, we have:
\begin{eqnarray}
\min_{u \succeq 0}\frac{1}{2\eta_0} \|w_0 - Tu\|^2 + \sum_{i=1}^{\tilde{n}}\frac{1}{2\eta_i}\|w_i - H^iu\|^2 
\end{eqnarray}
This is a non-negative least squares problem, which we solve using the Fast Non-Negative Least Squares algorithm \cite{bro1997fast}.
Algorithm \ref{algo:fixed_N_general} summarizes all updates of this block-coordinate descent algorithm. There is only one block (with respect to $w_0$) that does not necessarily have a unique minimum, and so Algorithm~\ref{algo:fixed_N_general} converges to a stationary point by the results of ~\cite{tseng2001convergence}.
\begin{algorithm} 
\caption{Fluence Map Optimization with Fixed Parameters and Fractions}\label{algo:fixed_N_general} 
\begin{algorithmic}[1]
\State \textbf{Input:} $u^{(0)}, \eta_0, \eta_1, \dots, \eta_{\tilde n}, N_1, \dots, N_M$
\Function{LowerLevelFixedParams}{$u^{(0)}, \eta_0, \eta_1, \dots, \eta_{\tilde n},  N_1, \dots, N_M, I_{\text{mean}}, \tilde I_{\text{max}}$}
\State \textbf{Initialize:} $k=0$
\State $\tilde \alpha := \left[\begin{smallmatrix}
-N_1\alpha_1\\
-N_2\alpha_2\\
\vdots\\
-N_M\alpha_M
\end{smallmatrix}\right],
B := \left[\begin{smallmatrix}
N_1 \mbox{diag}(\beta_1) & & & \\
 & N_2 \mbox{diag}(\beta_2) & & \\
& & \ddots & \\
& & & N_M \mbox{diag}(\beta_M)\\
\end{smallmatrix}\right]$
\For{$i = 1, \dots, \tilde n$}
\State $\tilde \gamma^i := \left[\begin{smallmatrix}
N_1\gamma^i_1\\
N_2\gamma^i_2\\
\vdots\\
N_M\gamma^i_M
\end{smallmatrix}\right],
D^i := \left[\begin{smallmatrix}
N_1 \mbox{diag}(\delta^i_1) & & & \\
 & N_2 \mbox{diag}(\delta^i_2) & & \\
& & \ddots & \\
& & & N_M \mbox{diag}(\delta^i_M)\\
\end{smallmatrix}\right]$
\EndFor
\While{\text{not converged}}
    \State $k \gets k+1$
    \State $w_0^{(k)} \gets \mbox{prox}_{-\eta_0f} (Tu - \eta_0\tilde\alpha)$
    \For{$i \in I_{\text{mean}}$} 
        \State $\Omega_i = \{w_i: (\tilde \gamma^i)^Tw_i+w_i^TD^iw_i \le C^i_{\text{mean}} \} $
        \State $w_i^{(k)} \gets \mbox{proj}_{\Omega_i}(H^i u)$
    \EndFor
    \For{$i \in  \tilde I_{\text{max}}$} 
        \State $\Omega_i = \{w_i: (\tilde \gamma^i)^Tw_i+w_i^TD^iw_i \le C^i_{\text{max}} \} $
        \State $w_i^{(k)} \gets \mbox{proj}_{\Omega_i}(H^i u)$
    \EndFor
    \State $u^{(k)} \gets \arg\min_{u \succeq 0}{\tilde\alpha}^Tw_0^{(k)} - f(w_0^{(k)}) + \frac{1}{2\eta_0} \|w_0^{(k)} - Tu\|^2 + \sum_{i=1}^{\tilde n} \frac{1}{2\eta_i}\|w_i^{(k)} - H^i u\|^2 $
\EndWhile
\State $u^{*} \gets u^{(k)}$
\State \Return $u^*$
\EndFunction
\State \textbf{Output:} $u^*$
\end{algorithmic}
\end{algorithm}

\subsubsection{Automatic Parameter Selection}
\label{sec:parameters}
The relaxed problem ($\text{P4}$) uses multiple parameters: $\eta_0$ and $\{\eta_i\}_{i=1}^{\tilde n}$. The auxiliary variable 
$w_0$ corresponds to the tumor and $w_i$ corresponds to the $i$-th OAR with $\tilde{n}$ number of OARs in the problem (including the introduced one-voxel maximum-dose-constraint OARs). Decreasing $\eta_0$ enforces the optimality of the solution (i.e. maximizes the tumor BE) because $Tu$ is forced to be closer to $w_0$, while decreasing $\eta_i$ increases the penalty corresponding to the $i$-th OAR constraint and thus enforces the feasibility of our solution making $H^iu$ closer to $w_i$. 
The existence of the proximal operator of the function $(-f)$ imposes an upper bound on $\eta_0$:
\begin{eqnarray}
\label{eqn:threshold}
\eta_0 \le \frac{1}{2\max_i{B_{ii}}}.
\end{eqnarray}
That is, the optimality penalty $\frac{1}{2\eta_0}$ should be big enough for the proximal operator to exist. 
This gives us an initial value for the parameter choice procedure. Combining these ideas, we develop an algorithm for the automatic selection of the parameters $\eta_i$:
\begin{enumerate}
\item Start $\eta_0$ from the threshold in Equation (\ref{eqn:threshold}). Initialize $\{\eta_i\}_{i=1}^{\tilde{n}}$ with the same value. Find the solution $u$.
\item Check if any OAR constraints are violated by $u$ found in Step 1.
\item Enforce feasibility: If there are any violated constraints, decrease $\eta_i$ by setting $\eta_i^+ = \Delta_{\eta}\eta_i$ with $\Delta_\eta < 1$ and solve for new $u^+$. Repeat steps 2-3 until all OAR constraints are satisfied.
\item Enforce optimality, that is, decrease $\eta_0$ by setting $\eta_0^+ = \Delta_{\eta}\eta_0$ and resolve for $u^+$ until $u^+$ fails to satisfy any constraint within the required tolerance for the OAR constraint. 
\end{enumerate}

Algorithm \ref{algo:auto_param_general} summarizes the lower-level optimization solution algorithms including the automated parameter selection with the fixed fractionation schedule $N_1, \dots, N_M$.

\begin{algorithm} 
\caption{Lower Level Optimization with Automated Parameter Selection (fixed fractions)}\label{algo:auto_param_general} 
\begin{algorithmic}[1]
\State \textbf{Input:} $u^{(0)}, N_1, \dots, N_M$
\Function{LowerLevelSolver}{$u^{(0)}, N_1, \dots, N_M$}
\State \textbf{Initialize:} $\eta_0 = \eta_1 = \dots = \frac{1}{2\max_i{B_{ii}}}$
\State $u \gets$ \textsc{LowerLevelFixedParams}($u^{(0)}, \eta_0, \eta_1, \dots, \eta_{\tilde n}$)
\While{there exists a violated constraint}
    \For{i in the violated constraints index set}
        \State $\eta_i \gets \eta_i\cdot\Delta_\eta$ \Comment{Decrease $\eta_i$, enforce feasibility}
        \State $u \gets$ \textsc{LowerLevelFixedParams}($u^{(0)}, \eta_0, \eta_1, \dots, \eta_{\tilde n}$)
    \EndFor
\EndWhile
\While{there is no violated constraints}
\State $\eta_0 \gets \eta_0\cdot\Delta_\eta$ \Comment{Decrease $\eta_0$, enforce optimality}
\State $u \gets$ \textsc{LowerLevelFixedParams}($u^{(0)}, \eta_0, \eta_1, \dots, \eta_{\tilde n}$)
\EndWhile
\State $u^* \gets u$
\State \Return $u^*$
\EndFunction
\State \textbf{Output:} $u^*$
\end{algorithmic}
\end{algorithm}

\section{Numerical Simulations}
\label{sec:numerical}
We apply the proposed framework with two different radiation modalities, $M_1$ and $M_2$, to the 2D phantom geometry shown in Figure \ref{fig:geometry}.  Since photons are currently the most widely used radiation modality in practice, we investigate the impact of combining photons ($M_1$) with a second modality ($M_2$), which has distinctive dosimetric characteristics as shown in Figure \ref{fig:isodose}. In Section \ref{subsec:tp}, 
we explain treatment planning and evaluation. 
Our framework is applicable to an arbitrary combination of different modalities with unique dosimetric and biological characteristics, but to build intuition for the proposal we investigate the impact of the difference between two modalities on the optimal BE in simple stages. 
In Section \ref{subsec:dosimetricpower}, we 
consider the simple scenario of combining $M_1$ with $M_2$ that has a dosimetric difference only, i.e. 
all radiobiological parameters between $M_1$ and $M_2$ 
are identical but the dose mapping matrices are different, i.e., $T_1 \neq T_2$ and $H^i_1 \neq H^i_2$.  
In Section \ref{subsec:biologicalpower},
we also add a radiobiological difference between $M_1$ and $M_2$. 
We vary the tumor's linear coefficients in the LQ model for $M_2$ ($\alpha_2$) and then vary the differential in the damage done by $M_2$ between the tumor and OARs ($r)$. Finally, in Section \ref{subsec:optimalfrac}, we present the performance of the optimal fractionation algorithm by comparing it to the true solution found using the brute-force technique, where all possible combinations of $(N_1, N_2)$ are individually used to find an optimal solution.

In all studies, we compute the initial guess, $u_0$ in Algorithm \ref{algo:optimize_N}, to give a uniform dose of 70 Gy to the tumor without any OAR constraints, which is commonly used in practice for head-and-neck tumors. The codes for our numerical simulations are available upon request.

\subsection{Phantom Geometry and Treatment Planning}
\label{subsec:tp}
\begin{figure}
    \centering
    \includegraphics[width=4in]{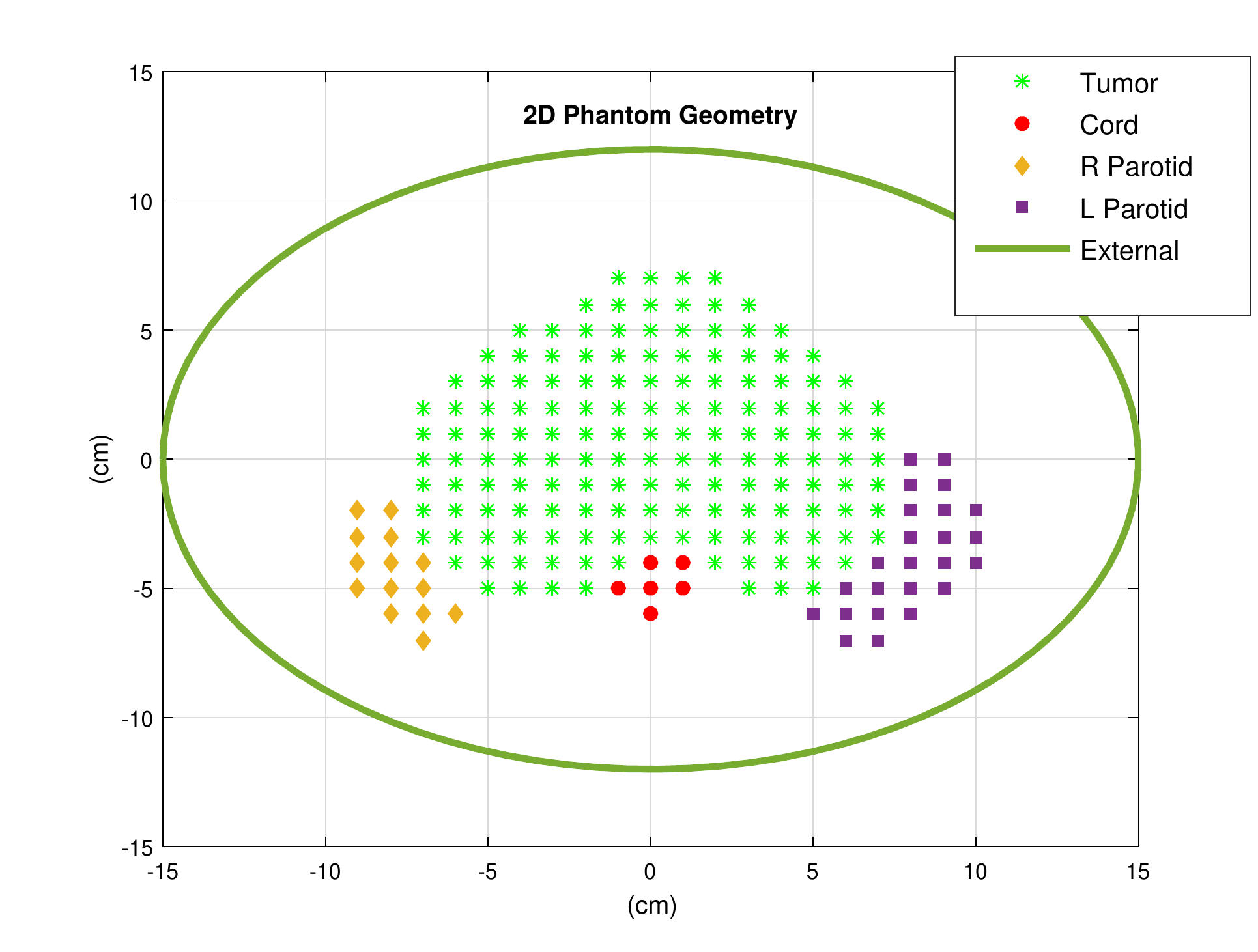}
    \caption{Phantom geometry}
    \label{fig:geometry}
\end{figure}

The 2D phantom geometry in Figure \ref{fig:geometry}
reflects a head-and-neck tumor surrounded by the spinal cord, right parotid, and left parotid glands.  
The unspecified tissue 1 cm inside of the external contour represents the skin of the patient. 
Adding a constraint for the unspecified tissue in the optimization ensures that the dose outside the tumor and OAR does not exceed the tolerance level. 

The dose mapping matrices for the first modality $(T_1, H^i_1$) were computed using the Elekta linear accelerator with 6 MV photons at the University of Washington.  
Radiation dose using $M_1$ is assumed to be delivered using seven equally spaced beams, 
i.e. gantry angles of 0$^{\circ}$, 51$^{\circ}$, 103$^{\circ}$, 153$^{\circ}$, 206$^{\circ}$, 257$^{\circ}$, and 309$^{\circ}$. 
The number of beamlets used is 195.  
The dose mapping matrices for the second modality ($T_2, H^i_2$) were computed using the proton beams with 250 MeV at the Seattle Cancer Care Alliance Proton Therapy Center with 40 spot positions within the tumor. 
The radiation dose using $M_2$ is assumed to be delivered in one beam (gantry angle of 0$^{\circ}$) 
as is often done in the proton therapy practice. 
We assume that the maximum number of fractions allowed is 25 fractions (i.e. $N_\text{max} = 25$).

Normal tissue tolerance BE and the radiobiological parameters used in computing BE were obtained from  published literature using photons and conventional fractionation schedules \cite{QUANTEC}. 
The $\alpha_1/\beta_1$ ratio is 10 Gy for the tumor with reference modality $M_1$. 
The $\gamma_1/\delta_1$ ratio is 2 Gy for the cord and unspecified tissue, 
and 5 Gy for the right and left parotid glands with $M_1$. 
The linear coefficients $\alpha_1$ and $\gamma_1$ of the LQ model are all assumed to be 0.35 Gy$^{-1}$ for $M_1$. 
The tolerance BE for all OARs used is summarized in Table \ref{tab:tolerance}.

\begin{table}[h!]
    \centering
    \begin{tabular}{rccc}
         & Constraint type & $\delta_1/\gamma_1$ & Tolerance BE \\
         \hline
     Cord    & Maximum dose & 2 Gy & 35\\
     Right Parotid     &  Mean dose & 5 Gy & 12\\
     Left Parotid     & Mean dose & 5 Gy & 12\\
     Unspecified tissue     & Maximum dose &  2 Gy & 13.125
    \end{tabular}
    \caption{Constraint type and tolerance biological effect for organs-at-risk (OAR) used in treatment planning.}
    \label{tab:tolerance}
\end{table}

We introduce a relative damage factor of $r = \delta_2/\alpha_2$ for $M_2$ to capture the relative effect of $M_2$ on the tumor and OAR. 
When a modality has a larger biological effect on the tumor, it also damages OAR more.  
However, the magnitude of the damage depends on tissue type, cell cycles, and other conditions. Therefore, $r$ is used to capture the differential in the damage between the tumor and OAR by $M_2$.  
As $r$ increases, the damage to OAR is relatively greater than the damage to the tumor, 
making $M_2$ clinically undesirable. We used $r \in \{0.5, 1.0, 1.5\}$.

The evaluation of the proposed framework is performed using percentage improvement in BE 
compared to BE obtained from two other methods with a single modality: 
(1) a conventional treatment course and (2) a single modality treatment course with an optimal fractionation schedule.  

Denote by $T_\text{conv}$ the conventional treatment course with fixed 25 fractions using photons only, i.e., $N_1 = 25$ fixed. 
Denote by $T_\text{single}$ the treatment course with a single modality with photons and an optimal number of fractions $N^\dag_1$.  
Then the percentage improvement of the objective function value ($=$ BE) using the optimal double-modality treatment course $T^*$ with $(N^*_1, N^*_2)$, where $1 \leq N^*_1 + N^*_2 \leq 25$, is given by
\begin{itemize}
\item Relative to a conventional course $T_\text{conv}$ with $N_1=25$ fixed
\[
\text{pObj}_\text{conv} = \frac{\mbox{BE using } T^*}{\mbox{ BE using } T_\text{conv}} \times 100 \%.
 \] 
 \item Relative to a single modality course $T_\text{single}$ with the optimal $N^*_1$
\[
\text{pObj}_\text{single} = \frac{\mbox{BE using } T^*}{\mbox{ BE using } T_\text{single}} \times 100 \%.
 \] 
\end{itemize}

\subsection{Dosimetric Difference}
\label{subsec:dosimetricpower}

In this section, we investigate the scenario where the only difference between $M_1$ and $M_2$ is the dose mapping matrices for the tumor and OAR.  
All other radiobiological parameters are set to be the same for the two modalities.  
The tumor doubling time, $T_d$, is varied between 2 and 100 days, that is, $T_d \in \{2,5,20,50,100\}$. The percentage improvement in the optimal BE with $T^*$ ranges between 105.7 \% and 107.9 \% compared to $T_\text{single}$, and between 107.7 \% and 122.9 \% compared to $T_\text{conv}$ depending on $T_d$.
The optimal number of fractions increases as $T_d$ increases, which agrees with clinical intuition. 
Slowly growing tumors (i.e. larger $T_d$) benefit from a long treatment course to reduce the long-term normal tissue side effect since normal tissues have better capability to recover from radiation damage than a tumor between fractions \cite{Hall2006}. 
The additional improvement seen in pObj$_\text{conv}$ compared to pObj$_\text{single}$ implies that the benefit of using $T^*$ comes from both using optimal fractionation schedule and multiple modalities. However, for a tumor with a large $T_d$, there is no benefit of using two modalities compared to a single modality with an optimal fractionation schedule. 
The complete results are shown in Table \ref{tab:protons}. 

\begin{table}[h!]
\begin{tabular}{m{1cm}m{3cm}m{3cm}m{2.5cm}m{2.5cm}}
   $T_d$  & Dual modality   & Single modality  &  &  \\
  (days) & optimal $(N^*_1, N^*_2)$ & optimal ($N^\dag_1$) & $\text{pObj}_\text{single}$ (\%)&  $\text{pObj}_\text{conv}$ (\%) \\
   \hline \hline \\
2 & (2, 2)& 2   & 105.7 & 122.9\\
\hline \\
5 & (6, 6)& 6  & 106.8 & 110.0 \\
\hline \\
10 & (13, 12) & 13  & 107.9 & 108.2\\
\hline \\
50 & (12, 12) & 25  & 107.8 & 107.8 \\
\hline \\
100 &(12, 12) & 25  & 107.8 & 107.8 \\
\hline \\
\end{tabular}
\caption{Optimal BE improvement with various tumor doubling time.  Parameters are fixed at $\alpha_2 = 0.35 \text{ Gy}^{-1}, r=1.0$}
\label{tab:protons}
\end{table}

\subsection{Radiobiological Difference}
\label{subsec:biologicalpower}
In this section we present the results when $M_2$ has a radiobiological difference from $M_1$ in addition to the dosimetric difference. 
First, we increase the linear coefficients in the LQ model, which means that $M_2$ kills more cells than $M_1$. 
We used $\alpha_2 \in \{0.35, 0.55, 0.75\}$ with $r$ fixed at 1. 
The percentage improvement in BE with $T^*$ ranges between 106.6 \% and 107.2\% compared to the single modality course with optimal fractionation $T_\text{single}$ and between 109.6 \% and 110.2 \% compared to the conventional treatment course $T_\text{conv}$ with 25 fractions. 
There is no significant difference in the results using different $\alpha_2$ values, 
and therefore the benefit of using $M_2$ is not clear. 
This is likely because killing more cells applies to both tumor and OAR. 
The complete results are shown in Table \ref{tab:alpha2}.

\begin{centering}
\begin{table}[ht!]
 \begin{tabular}{m{1cm}m{3cm}m{3cm}m{2.5cm}m{2.5cm}}
   $\alpha_2$   & Dual modality  & Single modality  &  &  \\
 $(Gy^{-1})$  & optimal $(N^*_1, N^*_2)$ & optimal ($N^\dag_1$) & $\text{pObj}_\text{single}$ (\%)&  $\text{pObj}_\text{conv}$ (\%) \\
   \hline \hline \\
0.35 & (6, 6)& 6  & 106.8 & 110.0 \\
\hline \\
0.55 & (10, 10)& 8  & 106.6  & 109.6 \\
\hline \\
0.75 & (8, 7) & 8 & 107.2  &  110.2\\
\hline \\
\end{tabular}
\caption{Optimal BE improvement with various $\alpha_2$: parameters are fixed at $r=1, T_d=5$ days.}
\label{tab:alpha2}
\end{table}
\end{centering}

We next varied the relative damage factor $r$, that is the differential in the damage between the tumor and OAR by $M_2$. 
We used $r \in \{0.5, 1.0, 1.5\}$. 
When $r=1$, the tumor's linear coefficient in the LQ model is the same as the OAR's linear coefficient, i.e. $\alpha_2 = \delta_2$. A smaller $r < 1$ implies that $M_2$ damages the tumor more than OAR ($\alpha_2 < \delta_2$) and $r > 1$ implies the reverse relation.  
The percentage improvement in BE with $T^*$ ranges between 102.2 \% and 125.7\% compared to the single modality course with optimal fractionation $T_\text{single}$ and between 105.1 \% and 129.4 \% compared to the conventional treatment course $T_\text{conv}$ with 25 fractions fixed. 
The results are shown in Table \ref{tab:r}. 
As $r$ decreases, the benefit of using a combination of $M_1$ and $M_2$ increases since 
$M_2$ damages the tumor more than OAR. This agrees with clinical intuition.

\begin{centering}
\begin{table}
 \begin{tabular}{m{1cm}m{3cm}m{3cm}m{2.5cm}m{2.5cm}}
   $r$   & Dual modality  & Single modality  &  &  \\
   & optimal $(N^*_1, N^*_2)$ & optimal ($N^*_1$) & $\text{pObj}_\text{single}$ (\%)&  $\text{pObj}_\text{conv}$ (\%) \\
   \hline \hline \\
0.5 & (2, 19)&  5 & 125.4 & 129.0 \\
\hline \\
1.0 & (6, 6)&  6 & 106.8  &  110.0 \\
\hline \\
1.5 & (4, 4) & 8 & 102.2  & 105.1 \\
\hline \\
\end{tabular} 
\caption{Optimal BE improvement with various $r$: parameters are fixed at $T_d=5$ days and $\alpha_2 = 0.35/\text{Gy}^{-1}$.}
\label{tab:r}
\end{table}
\end{centering}

\subsection{Uncertainty of $M_2$}
A modality with superior dosimetric or radiobiological characteristics often comes with a larger degree of uncertainty.  For example, protons deposit almost no dose beyond the Bragg peak but the uncertainty in the location of the Bragg peak makes them less desirable since depositing a large dose at a position slightly off target may result in an unacceptable dose to the tumor. 
Therefore, in practice, one strategy to mitigate  uncertainty is to ensure that a larger area around the tumor receives adequate dose by adding an extra margin around the clinical target volume.  
Delivering a larger dose in a bigger target volume often increases the dose to OAR, 
which is disadvantageous. 
To investigate the effect of this extra margin used for $M_2$ in our formulation, 
with hard constraints on OAR BEs, 
we add an extra margin to the tumor for $M_2$, compute the average BE per tumor voxel, and use it as an evaluation criteria. 
The average BE per tumor voxel for $M_1$ is computed without the extra margin, 
and therefore we expect $M_1$ to contribute a larger BE to the total tumor BE compared to the scenario, 
where the extra margin is not used for $M_2$.  
The results of using the average tumor BE per voxel to compute the pObj$_{\text{single}}$ and pObj$_{\text{conv}}$ reflect this clinical intuition and the optimal BE improvement of using multi-modality is less than what was achieved in Section \ref{subsec:biologicalpower}.  
The results are shown in Tables \ref{tab:ringTd}-\ref{tab:ringr}.

\begin{table}[h!]
\begin{tabular}{m{1cm}m{3cm}m{3cm}m{2.5cm}m{2.5cm}}
   $T_d$  & Dual modality   & Single modality  &  &  \\
  (days) & optimal $(N^*_1, N^*_2)$ & optimal ($N^\dag_1$) & $\text{pObj}_\text{single}$ (\%)&  $\text{pObj}_\text{conv}$ (\%) \\
   \hline \hline \\
2 &  (1, 1)&  2 &   102.5& 119.2 \\
\hline \\
5 &  (3, 4) &  3 & 103.1  & 105.7  \\
\hline \\
10 &  (12, 13) &    12&   102.9 &  103.2\\
\hline \\
50 &  (10, 14) &   25 &   103.1 &  103.1 \\
\hline \\
100 & (11, 14)&  25 &   103.4 &  103.4\\
\hline \\
\end{tabular}
\caption{Optimal BE improvement for a range of tumor doubling times when an extra margin to the tumor is used for $M_2$.  Parameters are fixed at $\alpha_2 = 0.35 \text{ Gy}^{-1}, r=1.0$}
\label{tab:ringTd}
\end{table}

\begin{centering}
\begin{table}[ht!]
 \begin{tabular}{m{1cm}m{3cm}m{3cm}m{2.5cm}m{2.5cm}}
   $\alpha_2$   & Dual modality  & Single modality  &  &  \\
 $(Gy^{-1})$  & optimal $(N^*_1, N^*_2)$ & optimal ($N^\dag_1$) & $\text{pObj}_\text{single}$ (\%)&  $\text{pObj}_\text{conv}$ (\%) \\
   \hline \hline \\
   0.35 &  (3, 4)&  3 &   103.2 & 105.7 \\
\hline \\
0.55 &  (6, 8) &  6 & 102.3  & 105.4  \\
\hline \\
0.75 &  (4, 5) &  4 &   103.2 &  106.0\\
\hline \\
\end{tabular}
\caption{Optimal BE improvement with various $\alpha_2$ when an extra margin to the tumor is used for $M_2$: parameters are fixed at $r=1, T_d=5$ days.}
\label{tab:ringalpha2}
\end{table}
\end{centering}

\begin{centering}
\begin{table}
 \begin{tabular}{m{1cm}m{3cm}m{3cm}m{2.5cm}m{2.5cm}}
   $r$   & Dual modality  & Single modality  &  &  \\
   & optimal $(N^*_1, N^*_2)$ & optimal ($N^*_1$) & $\text{pObj}_\text{single}$ (\%)&  $\text{pObj}_\text{conv}$ (\%) \\
   \hline \hline \\
      0.8 &  (7, 11)&  6 &   104.7 & 107.9 \\
\hline \\
1 &  (3, 4) &  3 & 103.2  & 105.7  \\
\hline \\
1.2 &  (3, 3) &  3 &   102.0 &  104.7\\
\hline \\
\end{tabular} 
\caption{Optimal BE improvement with various $r$ when an extra margin to the tumor is used for $M_2$: parameters are fixed at $T_d=5$ days and $\alpha_2 = 0.35/\text{Gy}^{-1}$.}
\label{tab:ringr}
\end{table}
\end{centering}

\subsection{Optimal Fractionation}
\label{subsec:optimalfrac}
\begin{figure}
\begin{centering}
\begin{tabular}{c}
\includegraphics[width=0.5\textwidth]{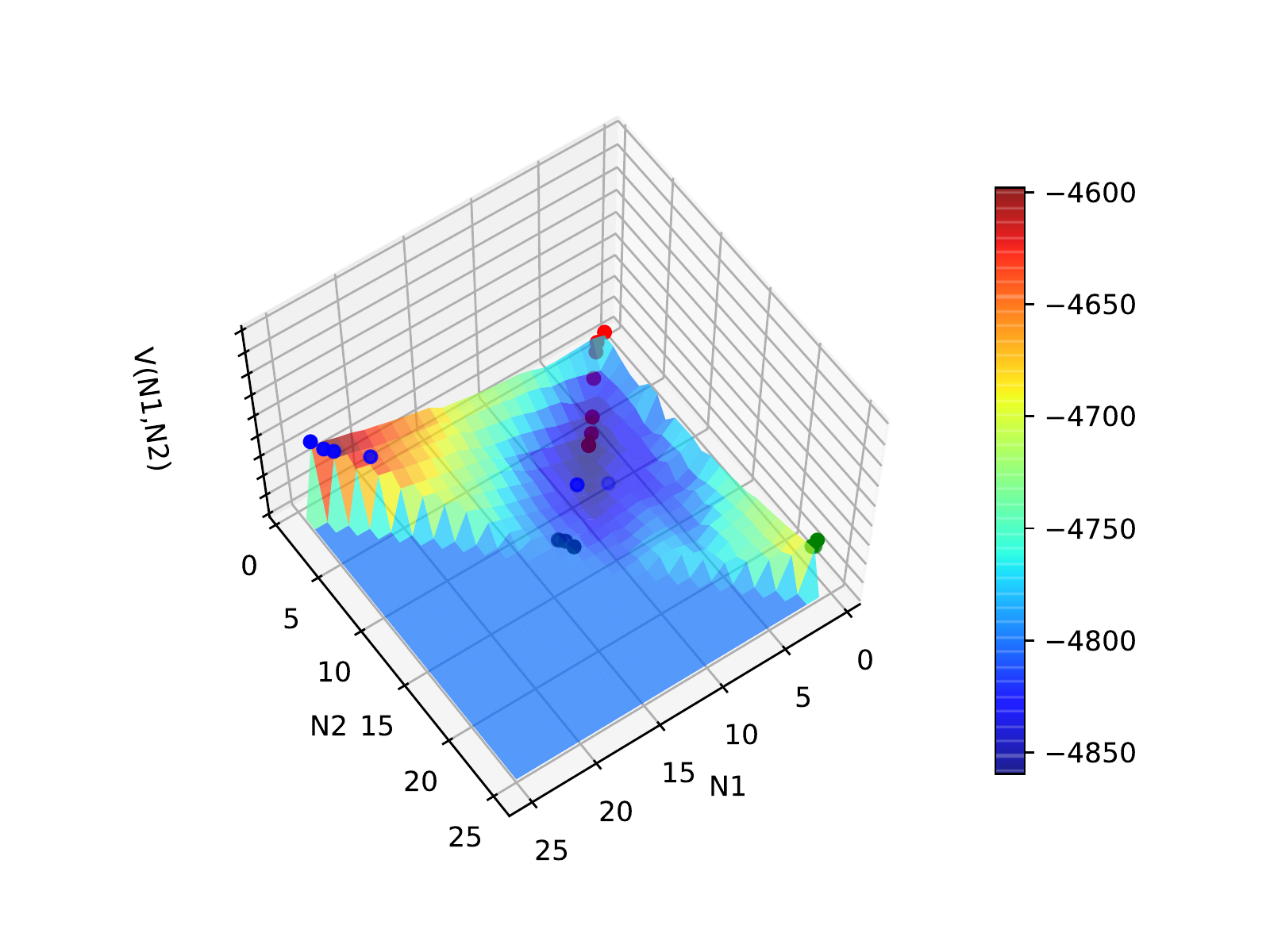}
\includegraphics[width=0.4\textwidth]{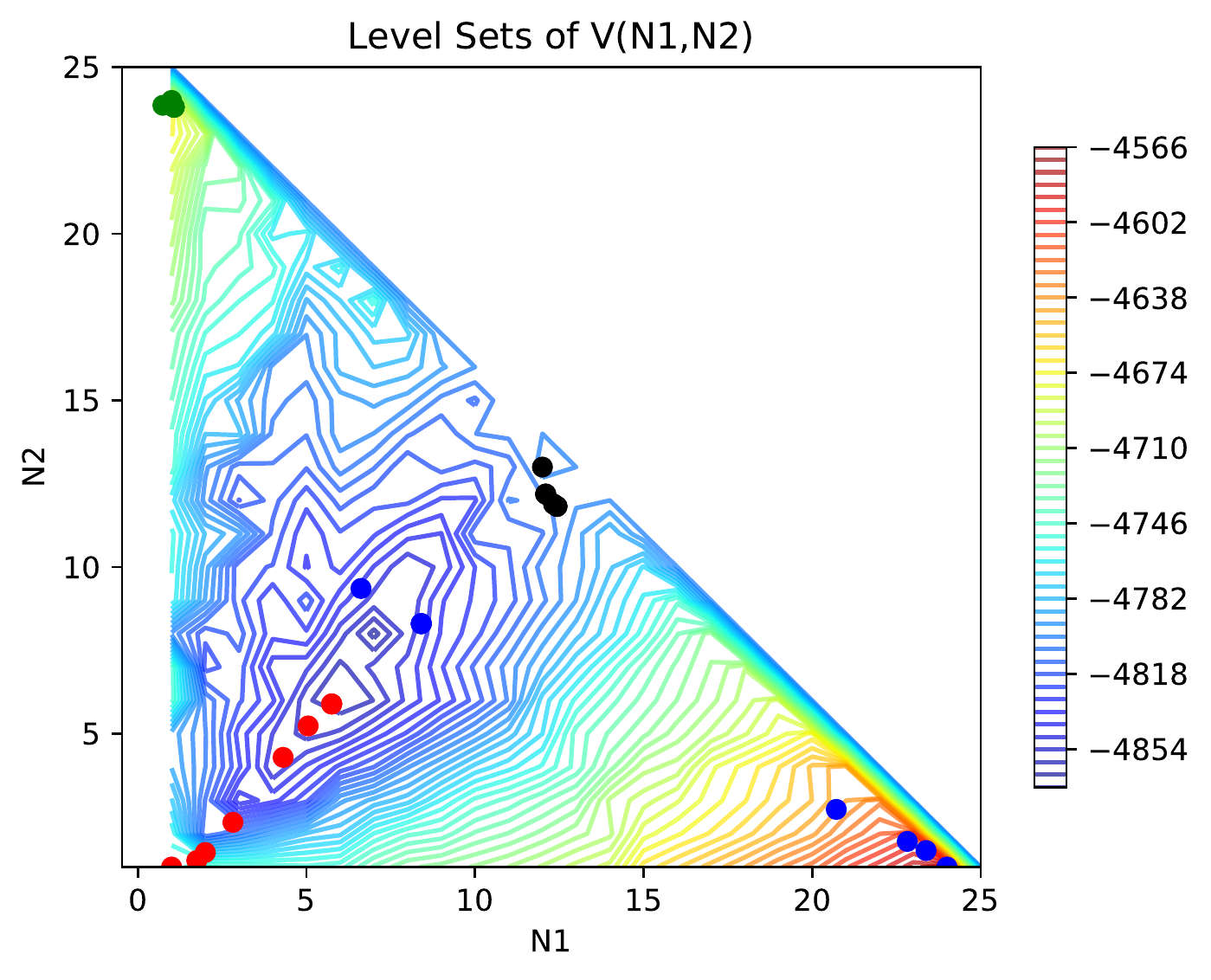}
\end{tabular}
\end{centering}
\caption{The surface (left) represents $V(N_1, N_2)$ computed by brute force. The dots represent iterates starting from four different initial guesses.}
\label{fig:optimizeN}
\end{figure}

In this section, we demonstrate the performance of our approach for the fractionation schedule optimization. We compute the ``ground truth" using the brute-force technique, where ($\text{P3}$) is solved for all possible integer ($N_1, N_2$) combinations with the constraint $1 \leq N_1+N_2 \leq N_\text{max}$ (that is, we sample the value function $V(N_1, N_2)$ for all the points of the integer grid $1 \leq N_1+N_2 \leq N_\text{max}$).
Figure \ref{fig:optimizeN} shows the value function $V(N_1, N_2)$ defined in \eqref{eqn:valuefunc} in 3D (left) and the level sets of $V(N_1, N_2)$ (right) obtained from the brute-force technique. The iterates $(N_1, N_2)$ of Algorithm \ref{algo:optimize_N} starting from the following four different initial guesses are shown as colored dots on Figure \ref{fig:optimizeN} (the colors correspond to the different initial points):
\begin{enumerate}
    \item ($N^0_1, N^0_2) = (1, N_\text{max}-1)$
    \item ($N^0_1, N^0_2) = (N_\text{max} -1, 1)$
    \item ($N^0_1, N^0_2) = (\lfloor N_\text{max}/2 \rfloor, \lceil N_\text{max}/2 \rceil)$
    \item ($N^0_1, N^0_2) = (1, 1)$
\end{enumerate}
As shown in Figure \ref{fig:optimizeN}, the output of Algorithm \ref{algo:optimize_N} is dependent on the initial guess, however, simple initial guesses as implemented in our algorithm may be sufficient to find a clinically relevant optimal solution. In our experiments, the above 4 initial guesses and 8 iterations of the upper level Algorithm \ref{algo:optimize_N} were sufficient to identify solutions close to the optimal.

\section{Conclusion}
\label{sec:conclusion}

We have developed a novel framework for radiation treatment planning with two or more radiation modalities. The principle idea is to maximize the biological effect using the multiple modalities by exploiting their unique dosimetric and biological characteristics captured through the dose mapping matrices and radiobiological parameters in the LQ dose-response model. Our framework allows $N_m=0$ as long as the sum of fractions of all modalities is equal to or greater than 1 fraction, and therefore it will correctly identify an optimal combination even when a single modality is optimal. The proposed framework ultimately offers an opportunity for an optimal multi-modality treatment planning paradigm. The feasibility of the proposed method was demonstrated using a simple phantom case and two modalities with varying only one parameter at a time such that clinical intuition can be applied. The agreement of the outcome with clinical intuition validates the potential use of our algorithm in more clinically relevant, complicated scenarios, where the clinical intuition is not readily available. 

The challenge of this approach is the dependence  of optimal solutions on the radiobiological parameters, which is indeed a challenge of biological treatment planning in general. Active research on advanced imaging could relieve some of the concerns in estimating these parameters \cite{Matuszak2019, Fiorentino2019}.  A novel algorithm we proposed to optimize the fractionation schedule with multiple modalities can be easily adapted for a single modality case.  A good choice of the initial guess is prudent for the optimal fractionation schedule algorithm to converge to an optimal solution. Applications of our algorithm to actual patient cases to evaluate the clinical significance and learning a good initial guess in the optimal fractionation schedule algorithm from a large patient database are left for future work.


\section*{Acknowledgment} This research was funded in part by the National Science Foundation through grant CMMI \# 1560476, and by the Washington Research Foundation Data Science Professorship. 

\begin{appendices}

\section{Proximal Operator Calculation}
\label{appendix:proximal}

\begin{eqnarray}
\label{eqn:prox}
\mbox{prox}_{-\eta_0f}(y) = \arg\min_x \left\{-f(x) + \frac{1}{2\eta_0}\|x-y\|^2 \right\}.
\end{eqnarray}

Equation (\ref{eqn:w0update}) gives us an update $w_0^+$ for $w_0$. We  note that the proximal operator is defined for convex functions whereas our $-f(x)$ is concave. Therefore, the minimization problem in (\ref{eqn:prox}) may be unbounded.  We will find the proximal operator by definition, that is, we will directly solve the minimization problem (\ref{eqn:prox}). Recall that $B$ is diagonal.
\begin{eqnarray}
&&\min_x -x^TBx + \frac{1}{2\eta_0}\|x-y\|^2 \\
&& \Leftrightarrow  
\min_x -x^TBx + \frac{1}{2\eta_0}x^Tx - \frac{1}{\eta_0}\langle y, x\rangle\\
&& \Leftrightarrow
\min_x -\sum_i B_{ii}x_i^2 + \frac{1}{2\eta_0}\sum_ix_i^2 - \frac{1}{\eta_0}\sum_i y_ix_i.
\end{eqnarray}
Finally, it is equivalent to
\begin{eqnarray}
\min_x -\sum_i (B_{ii}-\frac{1}{2\eta_0})x_i^2 - \frac{1}{\eta_0}\sum_i y_ix_i.
\end{eqnarray}
Consider possible cases:
\begin{itemize}
\item If $\exists i: (B_{ii}-\frac{1}{2\eta_0}) > 0, \text{then take } x: x_j=0 \; \forall j \not= i, x_i \rightarrow \infty$, the objective is unbounded below since the quadratic term dominates over the linear term.
\item If $\max_i(B_{ii}-\frac{1}{2\eta_0}) \le 0$, then we have a convex objective and therefore we can set derivative equal to zero to find the optimum. We have:
$$
-2Bx + \frac{1}{2\eta_0}2(y-x)(-1) = 0 \Rightarrow (-2B+\frac{1}{\eta_0}I)x = \frac{1}{\eta_0}y.
$$
Thus,
$$
(\mbox{prox}_{-\eta_0f}(y))_j = \begin{cases} 
\infty, & \max_i(B_{ii}-\frac{1}{2\eta_0}) > 0 \;\&\; j = \arg\max_i(B_{ii}-\frac{1}{2\eta_0})\\
0, & \max_i(B_{ii}-\frac{1}{2\eta_0}) > 0 \;\&\; j \not= \arg\max_i(B_{ii}-\frac{1}{2\eta_0})\\
\frac{1}{\eta_0}y_j/(-2B_{jj}+\frac{1}{\eta_0}), & \max_i(B_{ii}-\frac{1}{2\eta_0}) \le 0.
\end{cases}
$$
\end{itemize}

\section{Projection Calculation}
\label{appendix:projection}
To find the projection, let us reformulate the problem as follows:
\begin{eqnarray*}
&\min_{w}  \|w - v\|^2 \\
&\text{s.t.} \; \gamma_1  w_1 + \gamma_2  w_2+ \dots + \gamma_n w_n + D_1w_1^2 + D_2w_2^2 + \dots + D_nw_n^2 \le C_{\text{mean/max}}
\end{eqnarray*}
where $D_i$ stands for the $i$-th diagonal element of $D$. Now, let us complete the square:
\begin{eqnarray*}
&\gamma_1  w_1 + \gamma_2 w_2 + \dots + \gamma_n w_n + D_1w_1^2 + D_2w_2^2 + \dots + D_nw_n^2 \le C_{\text{mean/max}} \Leftrightarrow \\
&(\frac{\gamma_1}{2\sqrt{D_1}})^2 + 2(\sqrt{D_1}w_1)(\frac{\gamma_1}{2\sqrt{D_1}}) + (\sqrt{D_1}w_1)^2 + \dots \le C_{\text{mean/max}} +(\frac{\gamma_1}{2\sqrt{D_1}})^2 + \dots \Leftrightarrow \\
&\sum_i (\frac{\gamma_i}{2\sqrt{D_i}} + \sqrt{D_i}w_i)^2 \le C_{\text{mean/max}} + \sum_i(\frac{\gamma_i}{2\sqrt{D_i}})^2.
\end{eqnarray*}
Now, denote 
\begin{eqnarray*}
&K = C_{\text{mean/max}} + \sum_i(\frac{\gamma_i}{2\sqrt{D_i}})^2, \; z_i = \frac{\gamma_i}{2\sqrt{D_i}} + \sqrt{D_i}w_i, \text{ then }\\
&w_i = \frac{1}{\sqrt{D_i}}(z_i - \frac{\gamma_i}{2\sqrt{D_i}}).
\end{eqnarray*} 
We finally get:
\begin{eqnarray*}
&\min_z \sum_i(\frac{1}{\sqrt{D_i}}z_i - \frac{\gamma_i}{2D_i} - v_i)^2\\
&\text{s.t. } \|z\|^2 \le K
\end{eqnarray*}
Take the substitution $\hat D = \mbox{diag}(\frac{1}{\sqrt{D_1}}, \frac{1}{\sqrt{D_2}}, \dots, \frac{1}{\sqrt{D_n}}), \; l_i = \frac{\gamma_i}{2D_i} + v_i$ and consider the following (using KKT):
\begin{eqnarray*}
\min_z \frac{1}{2}\|\hat Dz - l\|_2^2 + \lambda(\|z\|^2 - K)
\end{eqnarray*}
In fact, we have:
\begin{eqnarray*}
L(z, \lambda) = \|\hat Dz - l\|_2^2 + \lambda(\|z\|^2 - K)\\
\text{KKT:}\\
\nabla L(z,\lambda) = \hat D^T(\hat Dz - l) + 2\lambda z = 0 \\
\lambda(z^Tz - K) = 0\\
z^Tz \le K \\
\lambda \le 0
\end{eqnarray*}
From the first condition, we have 
$$
z^* = (\hat D^T\hat D + \lambda I)^{-1}\hat D^Tl.
$$
If $\lambda = 0$, then we check the $z^*=(\hat D^T\hat D)^{-1}\hat D^Tl$ to satisfy ${z^*}^Tz^*\le K$. Otherwise, we look for the root $\lambda^*$ of
$$
\|(\hat D^T\hat D + \lambda I)^{-1}\hat D^Tl\|_2^2 = K,
$$
which we can find numerically by a root-finder routine. Then the projection is given by:
$$
z^* = (\hat D^T\hat D + \lambda^* I)^{-1}\hat D^Tl.
$$
Going back to the original variables,
$$
w_i = \frac{1}{\sqrt{D_i}}(z^*_i - \frac{\gamma_i}{2\sqrt{D_i}}).
$$

\end{appendices}

\bibliographystyle{unsrt}
\bibliography{Master}

\end{document}